\documentclass[a4paper,12pt]{article}
\usepackage{amsmath,amsfonts,amssymb,amsthm}

\usepackage{graphics}
\usepackage[dvips]{graphicx}

\newcommand{\ep}{\ensuremath{\varepsilon}}
\newcommand{\De}{\ensuremath{\Delta}}

\begin{document}

\begin{center}

\large

{\bf A dynamical model of opinion formation in voting \\ processes under bounded confidence}

\end{center}

Sergei Yu. Pilyugin\footnote{supported by the Russian Foundation for
Basic Research, grant 18-01-00230A},\\
St.Petersburg State University, Universitetskaya nab., 7-9, \\St.Petersburg, 199034, Russia. E-mail: {\sf sergeipil47@mail.ru}

and

M.C.~Campi\footnote{supported by the Russian Science Foundation,
grant 16-19-00057},\\
Dipartimento di Ingegneria dell'Informazione, Universit\`a di Brescia, \\via Branze 38, 25123 Brescia, Italia. E-mail: {\sf marco.campi@unibs.it}
\bigskip

{\bf Abstract.}
\noindent
In recent years, opinion dynamics has received an increasing attention, and various models have been introduced and evaluated mainly by simulation. In this study, we introduce and study a dynamical model inspired by the so-called ``bounded confidence'' approach where voters engaged in an electoral decision with two options are influenced by individuals sharing an opinion similar to their own. This model allows one to capture salient features of the evolution of opinions and results in final clusters of voters. 
The model is nonlinear and discontinuous. We provide a detailed study of the model, including a complete classification of fixed points of the appearing dynamical system and analysis of their stability. It is shown that any trajectory tends to a fixed point. The model highlights that the final electoral outcome depends on the level of interaction in the society, besides the initial opinion of each individual, so that a strongly interconnected society can reverse the electoral outcome as compared to a society with looser exchange. 
\bigskip

{\bf Introduction}
\bigskip

Studies on opinion dynamics aim to describe the processes by which opinions develop and take form in social systems, 
and research in this field goes back to the early fifties [8, 12]. In opinion studies, the word ``consensus'' refers to the agreement among individuals of a society towards a common view, a concept relevant to diverse endeavors of societal, commercial, and political interest. Consensus in opinion dynamics has been the object of several contributions such as [11, 24, 25, 29, 4, 5, 17]. A commonplace of these studies is that public opinion often evolves to a state in which one opinion predominates, but complete consensus is seldom achieved. Some basic models to describe opinion dynamics are described in the recent monographs [28] and [22].
\medskip

Most models in opinion dynamics are linear. One of the first nonlinear models was analyzed in [20, 19], where the notion of ``bounded confidence'' was also introduced. Bounded confidence concepts were further developed in [14, 6], while other nonlinear models based on similar approaches were studied in [9, 10]. In 2002, Hegselmann and  Krause [15] published an interesting study about an opinion model with bounded confidence, later called the Hegselmann -- Krause (HK) model, and provided computer simulations to illustrate the behavior of this model. In the same publication, they also noted that ``rigorous analytical results are difficult to obtain." After that, the HK model and its generalizations attracted a significant deal of attention, see, e.g.,
[7, 27, 21, 1, 23, 3, 2, 18, 30, 31, 13].  In particular, some theoretical results on sufficient conditions of convergence valid for a wide class of models of continuous opinion dynamics based on averaging (including the HK model and some models studied by Weisbuch and Deffuant) were obtained in [26]. The paper [16] extends the HK model 
by also including leaders and radical groups and derives various interesting behaviors resulting from this extension. 
\medskip

In this paper, we are especially interested in the dynamics of voters which have to choose between two alternatives. In this context, a natural assumption is that voters are more influenced by individuals sharing a similar opinion, which, when taken to its extreme, leads to models with bounded confidence. We discuss in more detail this aspect below after introducing the model. We contend that this situation leads to fixed points in the dynamics that correspond to the formation of opinion clusters. We study analytically these fixed points also analyze their stability properties. Although the present study refers to a simplified model, it is able to unveil and explain at a theoretical level fundamental features which have been observed in practice.
\medskip

While the model is described in detail in the next section, for explanation purposes we feel advisable to introduce here certain salient features of it. 
A population is formed by $N$ individuals, also called ``agents.'' The agents' opinion in regard of an electoral question with two options (identified by the numbers $-1$ and $1$) is described by $v_k \in [-1,1]$, $k = 1,\dots,N$, where a value close to $-1$ means that the individual $k$ carries an opinion more in favor of the option $-1$, while the opposite holds with a value $v_k$ close to $1$. Opinions $v_k$ evolve in discrete time through interaction. At any point in time, the new opinion of agent $k$ is formed by taking into account the opinions of agents whose values $v_l$ are not too distant from $v_k$ (bounded confidence). More precisely, fix a number $\epsilon > 0$ (not necessarily a small number) and denote by $J(v_k)$ the set of indices $l$ of
agents with opinions $\ep$-close to $v_k$, i.e.
$$
J(v_k)=\{l\in\{1,\dots,N\}:\;|v_l-v_k|\leq\ep\}.
$$
The new opinion of agent $k$ is obtained by adding to $v_k$ a value proportional to the average of opinions $v_l$ over the set $J(v_k)$ and ``cutting'' the new value if it exceeds the boundaries of the interval $[-1,1]$ (for a precise description, refer to the next section). It turns out that, apart from special configurations that give unstable equilibria, this dynamics leads to final configurations where the population splits  in two clusters, having values $-1$ and $1$. When taking the average to compute the value by which $v_k$ is updated,  the value of agent $k$ is included in the calculation as well. In the extreme case where an agent has no other $\ep$-close agents, this implies that this agent reinforces her/his belief: in absence of counter-arguments, one tends to strengthen her/his own initial opinion; in general, one's opinion is compared with the opinion of others in a neighborhood to determine the evolution.
\medskip

In the proposed model, an agent is only influenced by agents who are having a similar idea. This modeling assumption only holds in first approximation as agents may also interact with others that think quite differently and get influenced by them. Hence, this model only captures the predominant elements in a social interaction, while it neglects various second-order aspects. We also note that assuming that agents are ``deaf'' to others thinking differently is getting more realistic as the world evolves towards interaction schemes based on social media and the web where the contacts and sources of information are selected by the users.
\medskip

The structure of the paper is as follows. In Section~1, the mathematical definition 
of the model is given. Section~2 is devoted to the behavior
of the dynamical system generated by the model: we describe fixed points, study their stability, and show that any positive trajectory tends to a fixed point. Numerical examples are finally presented in Section~3. These examples show interesting features, 
for example, that the level of interaction influences the opinions in the long run to the point that the predominance of one option over the other can be reverted depending on the interaction level in the society.
\bigskip

{\bf 1. Definition of the opinion model}
\bigskip

The opinion of $N$ agents is described by a finite array
$$
V = (v_k\in[-1,1]: \;k=1,\dots, N),
$$
where $v_k$ has to be interpreted as the level of appreciation of agent $k$ for one among two options: a value $v_k$ close to $-1$ means that agent $k$ has a preference for option $-1$, and the closer $v_k$ to $-1$, the stronger the preference; the opposite holds for option $1$. Denote by ${\cal V}=[-1,1]^N$ the set of such arrays.
\medskip

We fix two numbers $h,\ep\in(0,1)$. In addition, we fix two functions, $a(v,w)$ and $i(v)$ (called {\it affinity} and {\it influence} function, respectively).
\medskip

In this study, the function $a(\cdot,\cdot)$ is defined as follows:
\begin{equation}
\label{a}
\nonumber
a(v,w)=1\mbox{ if }|v-w|\leq \ep \mbox { and }a(v,w)=0 \mbox{ otherwise}.
\end{equation}
If $a(v_k,v_l)=1$, we sometimes say that ``$v_k$ is influenced by $v_l$.''

The function $i(\cdot)$ is defined simply as follows:
\begin{equation}
\label{i1}
i(v)=v.
\end{equation}

For $k=1,\dots,N$, denote by $J(v_k)$ the set of indices $l \in [1,N]$ (here and below, we denote by $[a,b]$ the set of indices $\{a,\dots,b\}$) such that $|v_l-v_k|\leq\ep$ and by $I(v_k)$ the cardinality
of the set $J(v_k)$. 
\medskip

We study the dynamics on ${\cal V}$ defined by the following operator $\Phi$. First we fix a $V\in{\cal V}$ and consider the auxiliary array
$$
W(V)=(w_1(V),\dots,w_N(V))
$$
defined as follows
\begin{equation}
\label{01}
\nonumber
w_k(V)=v_k+h\sum_{l=1}^N \frac{i(v_l)a(v_k,v_l)}{I(v_k)},\quad
k=1,\dots, N.
\end{equation}
Sometimes, when this does not lead to confusion, we write $W(V)=(w_1,\dots,w_N)$ instead of $W(V)=(w_1(V),\dots,w_N(V))$.

Due to (\ref{i1}),
\begin{equation}
\label{010}
w_k(V)=v_k+\frac{h}{I(v_k)}\sum_{l\in J(v_k)}v_l.
\end{equation}

After that, we define
$$
\Phi(V)=(v_1',\dots,v_N')
$$
by ``cutting'' the elements of $W(V)$ according to the rule
$$
v'_k=-1\mbox{ if }w_k<-1,\quad v'_k=1\mbox{ if }w_k> 1,
$$
and
$$
v'_k=w_k\mbox{ if }|w_k|\leq 1.
$$

Obviously,
\begin{equation}
\label{001}
\nonumber
\Phi({\cal V}) \subseteq {\cal V}.
\end{equation}

Note that if we replace in \eqref{010} $v_l$ by $v_l-v_k$ and take $h=1$, then we get the HK model.

Our main goal is to study fixed points of the operator $\Phi$ and their stability.
\bigskip

{\bf 2. Dynamics of the opinion model}
\bigskip

We start with an initial array $V=V^0$ with the following property:
\begin{equation}
\label{1.3}
\nonumber
 v^0_1\leq\dots\leq v^0_N.
\end{equation}

Of course, we have a lot of ways to ``numerate'' members of our group (for example, alphabetically). Since the set $J(v_k)$ depends not on indices of elements of $V$ but on their values, it follows from formula (\ref{010}) that $\Phi(V)$ is, in a sense, independent of numeration.
In our case, an ordering reflecting the initial preferences of the agents seems to be the most convenient.

Let
\begin{equation}
\label{3.1}
\nonumber
V^n=\Phi^n(V^0)=(v^{n}_1,\dots, v^{n}_N).
\end{equation}

First let us note some important properties of the operator $\Phi$.

We need a simple technical statement (for its proof, see, for example, item (i) of Lemma~2 in [19]).
\medskip

{\bf Lemma 1. }
{\em If
$$
x_1\leq\dots\leq x_n\leq y_1\leq\dots\leq y_m,
$$
then}
$$
\frac{x_1+\dots+x_n}{n}\leq
\frac{x_1+\dots+x_n+y_1+\dots+y_m}{n+m}\leq
\frac{y_1+\dots+y_m}{m}.
$$
\medskip

Take a sequence $V=(v_1,\dots,v_N)$ such that
$$
v_1\leq\dots\leq v_N
$$
and consider the ``increments''
$$
\De_k=w_k(V)-v_k.
$$

{\bf Lemma 2. }{\em The following inequalities hold}:
\begin{equation}
\label{incr}
\De_{k+1}\geq \De_k,\quad k=1,\dots,N.
\end{equation}
\medskip

{\em Proof.} Let $J(v_k)= [a,a+l]$ with $v_a\leq\dots\leq v_{a+l}$ and $J(v_{k+1})=[b,b+m]$ with $v_b\leq\dots\leq v_{b+m}$. By formula (\ref{010}),
$$
\De_k=h\frac{v_a+\dots+v_{a+l}}{l+1}
$$
and
$$
\De_{k+1}=h\frac{v_b+\dots+v_{b+m}}{m+1}.
$$

If $J(v_k)\cap J(v_{k+1})=\emptyset$ (which is equivalent to
the inequality $a+l<b$), then (\ref{incr}) obviously holds.

Otherwise, let $J(v_k)\cap J(v_{k+1})=[b,a+l]$;
it follows from Lemma~1 that
$$
\frac{v_a+\dots+v_{a+l}}{l+1}\leq
\frac{v_{b}+\dots+v_{a+l}}{a+l-b+1}\leq
\frac{v_b+\dots+v_{b+m}}{m+1},
$$
which completes the proof. \qed \\

The following statements are more or less obvious but since we use them
 many times, we formulate them separately.

Applying induction on $n$ based on Lemma~2, the following properties of $V^n=\Phi^n(V^0)$ can easily be established.
\medskip

{\bf Corollary~1.}
\begin{itemize}
\item[(a)] {\em Every array $V^n$ is nondecreasing};
\item[(b)] {\em If $v^n_k=1$, then $v^n_l=1$ for} $l>k$;
\item[(c)] {\em If $v^n_k=1$, then $v^m_k=1$ for} $m>n$.
\end{itemize}
\medskip

We do not explicitly formulate obvious analogs of items (b) and (c) for $v^n_k=-1$.
\medskip

Let us explain a step of the induction in proving item (a) when we pass from $n=0$ to $n=1$ (the other steps are similar). Inequalities \eqref{incr} for $\De_k=w_k(V^0)-v_k^0$ and the nondecreasing property of $V^0$ imply that the  $W(V^0)$ is 
nondecreasing; hence, $V^1$ is nondecreasing as well.

(b) follows from (a).

(c) If $v^0_k=1$, then $v^0_k$ is not influenced by negative $v^0_l$ (since $\ep<1$); hence, $w^1_kgeq 1$ and $v^1_k = 1$.

\medskip
We next move to the study of fixed points of $\Phi$. 
Recall that $P\in{\cal V}$ is a fixed point of $\Phi$ if $\Phi(P)=P$.

First, we mention a class of fixed points which is important for us (as we show below, almost all positive trajectories of $\Phi$ tend to such fixed points). Let $P=(-1,-1,\dots,-1,1,\dots,1)$, where the first $L$ entries equal $-1$ while the remaining equal $1$.
We do not exclude the cases of $P=(-1,\dots,-1)$ (in which $L=N$) and $P=(1,\dots,1)$ (in which we formally set $L=0$). Any such $P$ is a fixed point of $\Phi$. This follows from item (c) of Corollary~1 (and its analog for $v^n_k=-1$).

Let us call any such $P=(-1,-1,\dots,-1,1,\dots,1)$ a basic fixed point of $\Phi$. We are going to show that any basic fixed point is asymptotically stable for $\Phi$ (see Theorem~1).

Let us start with a simple statement which we often use below.
\medskip

{\bf Lemma~3.}
{\em If $v_k^0\geq\ep$, then there exists an $n_0\geq 0$ such 
that $v_k^n=1$ for} $n\geq n_0$.
\medskip

{\em Proof. } The condition $v_k^0\geq\ep$ implies 
that $v_k^0$ is not influenced by negative $v_l^0$. On the other hand, 
$k\in J(v_k)$, so that
$$
w_k(V^0)\geq v_k^0+\frac{h}{N}v_k^0\geq\ep\left(1+\frac{h}{N}\right).
$$

If $w_k(V^0)\geq 1$, then $v_k^1=1$, and our statement follows from item (c) of Corollary~1. Otherwise,
$$
w_k(V^1)\geq \ep\left(1+\frac{h}{N}\right)+
\frac{h\ep}{N}\left(1+\frac{h}{N}\right)>
 \ep\left(1+\frac{2h}{N}\right),
$$
and so on, which obviously implies our statement. \qed \\

The same reasoning shows that if $v_k^0\leq -\ep$, then there exists an $n_0\geq 0$ such that $v_k^n=-1$ for $n\geq n_0$.
\medskip

Introduce the following metric on ${\cal V}$: if
$$
V=(v_1,\dots,v_N)\mbox{ and }V'=(v_1',\dots,v_N'),
$$
set
$$
\rho(V,V')=\max_{1 \leq k\leq N}|v_k-v'_k|.
$$

{\bf Theorem~1. }
{\em Let $P$ be a basic fixed point. If
\begin{equation}
\label{1.10}
\rho(V^0,P)\leq 1-\ep,
\end{equation}
then there exists an $n_0$ such that}
\begin{equation}
\label{1.11}
\Phi^n(V^0)=P \quad for \ n\geq n_0.
\end{equation}
\medskip

{\em Proof. } Let $V^0=(v^0_1,\dots,v^0_N)$ satisfy inequality (\ref{1.10}). Then
$$
|v^0_k|\geq\ep,\quad k=1,\dots,N,
$$
and our theorem follows from Lemma~3
since the number of components of $V^0$ is finite. \qed
\medskip

{\bf Remark~1. } One can establish the convergence to basic fixed 
points under weaker conditions than (\ref{1.10}). Assume, for example,
that
$$
v_1^0\leq\dots\leq v_L^0<0<v_{L+1}^0\leq\dots\leq v_N^0
$$
and
$$
v_{L+1}^0-v_L^0>\ep.
$$
Then the same reasoning as in Lemma~3 shows that $\Phi^n(V^0)=P$ 
for some finite $n$, where $P$ is a basic point.
\medskip

There exist fixed points that are not basic; we show below that they are unstable. A simple example of such a fixed point is as follows.
Let $N=3$; clearly, $P=(p_1,p_2,p_3)=(-1,0,1)$ is a fixed point of $\Phi$. We first describe all possible nonbasic fixed points of $\Phi$.
\medskip

{\bf Theorem~2. }
{\em If $P$ is a nonbasic fixed point of $\Phi$, then either
\begin{equation}
\label{nb1}
P=(-1,\dots,-1,0,\dots,0,1,\dots, 1)
\end{equation}
or
\begin{equation}
\label{nb2}
P=(-1,\dots,-1,p_a,\dots,p_l,0,\dots,0,p_b,\dots,p_m,1,\dots, 1),
\end{equation}
where
\begin{equation}
\label{nb3}
-\ep<p_k<0,\quad k\in[a,l],
\end{equation}
\begin{equation}
\label{nb4}
0<p_k<\ep,\quad k\in[b,m],
\end{equation}
\begin{equation}
\label{nb5}
J(p_k)=[a,m],\quad k\in[a,m],
\end{equation}
and}
\begin{equation}
\label{nb6}
p_a+\dots+p_m=0.
\end{equation}
\medskip

{\em Proof. } It is clear that if $P$ is a nonbasic fixed point that does not have form (\ref{nb1}), then it has form (\ref{nb2}) with $p_a,\dots,p_l\in(-1,0)$ and $p_b,\dots,p_m\in(0,1)$.

Inequalities (\ref{nb3}) and (\ref{nb4}) follow from Lemma~3.

Let us prove the remaining statements.

Since $p_a<0$ and $P$ is a fixed point, $p_a$ is influenced by positive $p_i$ and it cannot be influenced by $p_i=1$. Hence, there exists an index $r(a)\in[b,m]$ such that either
\begin{equation}
\label{1.2.1}
J(p_a)=[1,r(a)]
\end{equation}
or
\begin{equation}
\label{1.2.2}
J(p_a)=[a,r(a)].
\end{equation}

Note that these cases are different only if $a>1$.

Since $P$ is a fixed point,
\begin{equation}
\label{1.2.3}
-(a-1)+p_a+\dots+p_{r(a)}=0
\end{equation}
in the first case and
\begin{equation}
\label{1.2.4}
p_a+\dots+p_{r(a)}=0
\end{equation}
in the second case.

It follows from (\ref{nb3}) that any $p_k$ with
$k\in[a,l]$ is influenced by $p_a$.

Thus, there exists an index $r(a+1)\in[b,m]$
such that either
\begin{equation}
\label{1.2.5}
J(p_{a+1})=[1,r(a+1)]
\end{equation}
or
\begin{equation}
\label{1.2.6}
J(p_{a+1})=[a,r(a+1)].
\end{equation}

We claim that

$\bullet$ if $a>1$, then (\ref{1.2.1}) implies (\ref{1.2.5});

$\bullet$ (\ref{1.2.2}) implies (\ref{1.2.6});

$\bullet$ in both cases (\ref{1.2.5}) and (\ref{1.2.6}),
$r(a+1)=r(a)$.

To prove the first claim, we note that if $a>1$ and (\ref{1.2.1}) holds, then
$$
p_a+\dots+p_{r(a)}=a-1>0,
$$
while if (\ref{1.2.6}) holds, then
$$
p_a+\dots+p_{r(a)}=0\mbox{ if }r(a)=r(a+1)
$$
and
$$
p_a+\dots+p_{r(a)}=-p_{r(a)+1} - \cdots -p_{r(a+1)}<0\mbox{ if }r(a)\neq r(a+1).
$$

The second claim follows from the fact that if $p_a$ is not
influenced by $p_i=-1$, then $p_{a+1}\geq p_a$ cannot be
influenced by $p_i=-1$ as well.

To prove the third claim, we compare the equality
$$
-(a-1)+p_a+\dots+p_{r(a)}+p_{r(a+1)}=0
$$
with (\ref{1.2.3}) in the first case and the equality
$$
p_a+\dots+p_{r(a)}+p_{r(a+1)}=0
$$
with (\ref{1.2.4}) in the second case and note that $p_{r(a+1)}$ must be positive.

Continuing this process, we conclude that either $J(p_k)=[1,r(k)]$ for all $k\in[a,l]$ or $J(p_k)=[a,r(k)]$ for all $k\in[a,l]$, and, in both cases,
$$
r(a)=r(a+1)=\dots=r(l).
$$

Clearly, this common value must be equal to $m$ (since, otherwise, $p_m$ is not influenced by negative $p_i$, which is impossible for the fixed point $P$).

In the second case, the equality $r(a)=m$ implies (\ref{nb5}), and equality (\ref{1.2.4}) implies (\ref{nb6}).

To complete the proof of the theorem, it remains to show that if $a>1$, then the first case is impossible.

To do this, let us start with $p_m$ and ``move in the opposite direction": find $t(m)\in[a,l]$ such that $J(p_m)=[t(m),N]$ or $J(p_m)=[t(m),m]$, and so on.

Repeating the above reasoning, we get either the equality
$$
p_a+\dots+p_m=m-N\leq 0
$$
or equality (\ref{nb6}); both contradict the equality
$$
p_a+\dots+p_m=a-1>0
$$
obtained above. \qed \\

Now we are going to prove that if $P$ is a nonbasic fixed point of $\Phi$, then $P$ is unstable in a strong sense: $P$ has a neighborhood $U$ such that for any point $V\in U$ not belonging to a subset of $U$ of positive codimension, the trajectory $\Phi^n(V)$ leaves $U$ as $n$ grows. The authors are grateful to A. Proskurnikov who have noticed this fact and suggested the idea of the proof of the following theorem.
\medskip

{\bf Theorem~3. }
{\em If $P$ is a nonbasic fixed point of $\Phi$ having form 
(\ref{nb1}) or (\ref{nb2}), then there exists a $d>0$ such that if
$$
U=\{V:\;\rho(V,P)<d\},
$$
and
\begin{equation}
\label{1.3.0}
\nonumber
\Pi=\{V:\;v_a+\dots+v_m=0\},
\end{equation}
then for any point $V\in U\setminus \Pi$ there exists an $n>0$ such that} $\Phi^n(V)\notin U$.
\medskip

{\em Proof. } We impose several conditions on $d$.

First, it follows from (\ref{nb3}) and (\ref{nb4}) that we can 
take $d$ so small that if $V\in U$, then
\begin{equation}
\label{1.3.1}
-\ep<v_a\leq\dots\leq v_m<\ep.
\end{equation}

Second, condition (\ref{nb5}) implies that if $a\neq 1$ 
(i.e., $P$ has components equal to $-1$), then $p_a$ is not influenced by these components (i.e., $p_a+1>\ep$). Similarly, $p_m$ is not influenced by 
components $+1$ (if they exist). Hence, we can take $d$ so small that 
if $V\in U$, then
\begin{equation}
\label{1.3.2}
\nonumber
J(v_k) \subset [a,m],\quad k\in[a,m].
\end{equation}

Finally, we take $d$ so small that
\begin{equation}
\label{1.3.3}
\frac{h}{N}(p_m-d)>2d \ \mbox{ and } \ -\frac{h}{N}(p_a+d)>2d
\end{equation}
(recall that $p_m>0$ and $p_a<0$).

Denote
$$
s(V)=v_a+\dots+v_m.
$$

First we claim that if $V\in U$ and
\begin{equation}
\label{1.3.4}
J(v_m)\neq [a,m],
\end{equation}
then $\Phi(V)\notin U$.

Assume that $s(V)\geq 0$. It follows from (\ref{1.3.1}) and (\ref{1.3.4}) that $J(v_m)=[k,m]$, where $k\leq b$ (since $v_m$ is influenced 
by all positive components of $V$ with indices in $[b,m]$).

If $w_m(V)\geq 1$, then $v^1_m$ equals $1$, and our claim follows 
from (\ref{1.3.1}). Otherwise,
$$
v^1_m=w_m(V)=v_m+\frac{h}{m-k+1}(v_k+\dots+v_m).
$$

Since
$$
v_k+\dots+v_m=s(V)-(v_a+\dots+v_{k-1})\geq -v_a
$$
(we take into account that $s(V)\geq 0$ and $v_{a+1},\dots,v_{k-1}\leq 0$), 
it follows from the inequalities $m-k+1\leq N$, $v_a<p_a+d$, and (\ref{1.3.3}) that
$$
v^1_m-v_m\geq -\frac{h}{m-k+1}v_a>-\frac{h}{N}(p_a+d)>2d,
$$
which is impossible if $V \in U$ and $\Phi(V)\in U$.

If $s(V)<0$, we apply a similar reasoning taking into account that relation (\ref{1.3.4}) implies the relation $J(v_a)\neq[a,m]$.

Now let us take a point $V\in U$ and assume that $\Phi^n(V)\in U$ for all $n>0$. It follows from our previous reasoning that in this case,
$$
J(v^n_k)=[a,m],\quad k\in[a,m],
$$
for all $n$. Then, denoting $s_0 = S(v)$, we get the equalities
$$
v^1_k=v_k+\frac{h}{m-a+1}(v_a+\dots+v_m)=
v_k+\frac{hs_0}{m-a+1},\quad k\in[a,m],
$$
which yields the equality
$$
s(\Phi(V))=s_0(1+h).
$$

Similarly,
$$
s(\Phi^2(V))=s(\Phi(V))(1+h)=s_0(1+h)^2,\dots , s(\Phi^n(V))=s_0(1+h)^n,
$$
and so on.

If $V\notin \Pi$, then $s_0\neq 0$, and the above value
is unbounded as $n\to\infty$,
which is impossible since the values
$S(V)$ for $V\in U$ are bounded.

This completes the proof. \qed \\

Now we prove that if
\begin{equation}
\label{1.4.0}
\ep\leq 1/2,
\end{equation}
then trajectories of $\Phi$ tend to fixed points.
\medskip

{\bf Theorem~4. }
{\em If condition (\ref{1.4.0}) is satisfied, then any trajectory 
$\Phi^n(V^0)$ tends to a fixed point of $\Phi$ as} $n\to\infty$.
\medskip

{\em Proof. }
Consider an initial sequence
$$
V=(v_1,\dots,v_N).
$$
Corollary~1 implies that if $v^n_k=-1$, then $v^m_l=-1$ for $m\geq n$
and $l\leq k$; similarly, if $v^n_k=1$, then $v^m_l=1$ for $m\geq n$
and $l\geq k$.

Since the number of components of $V^n = \Phi^n(V^0)$ is finite, we conclude that there exist integers $0\leq a<b\leq N$ and $n_1\geq 0$ such that 
if $n\geq n_1$, then
$$
V^n=(-1,\dots,-1,v^n_a,\dots,v^n_b,1,\dots,1),
$$
where
$$
|v^n_k|<1,\quad  k\in[a,b];
$$
in words, the number of components equal to $\pm 1$ ``stabilizes.''

If the ``middle'' part $(v^{n_1}_a,\dots,v^{n_1}_b)$ is absent, $\Phi^{n_1}(V)$ is a fixed point, and our statement is proved.

To simplify the notation, assume that $n_1=0$. Clearly, our problem is to describe the behavior of $(v^n_a,\dots,v^n_b)$ as $n$ grows.

It was shown in Lemma~3 that if $|v^0_k|\geq\ep$ for 
some $k\in[a,b]$, then $|v^n_k|=1$ for large $n$, which is impossible. Hence,
\begin{equation}
\label{1.4.1}
-\ep< v^n_a\leq\dots\leq v^n_b< \ep,\quad n\geq 0.
\end{equation}

These inequalities and  condition (\ref{1.4.0}) imply that
$$
J(v^n_k)\subset [a,b],\quad k\in[a,b],\;n\geq 0.
$$

It follows that the behavior of $(v^n_a,\dots,v^n_b)$ is determined by 
components of $V^n$ with indices from $a$ to $b$. Thus, without loss of generality, 
we may assume that we study the behavior of $V^n$ 
with $|v^n_k|<\ep,\;1\leq k\leq N, \quad n \geq 0$.

Let
$$
{\cal N}(V)=\{(k,l)\in[1,N]\times [1,N]:\;|v_k-v_l|>\ep\}
$$
be the set of pairs $(k,l)$ of indices such that $v_k$ is not influenced by $v_l$ et vice versa.

We prove the following simple but relevant statement separately.
\medskip

{\bf Lemma~4. }
\begin{equation}
\label{1.4.2}
{\cal N}(V^n)\subset {\cal N}(V^{n+1}),\quad n\geq 0.
\end{equation}
\medskip

{\em Proof. } Inclusion (\ref{1.4.2}) means that if
$|v^n_k-v^n_l|>\ep$, then
$$
|v^{n+1}_k-v^{n+1}_l|>\ep
$$
as well.

Assume, for the sake of clarity, that $v^n_k-v^n_l>\ep$ 
(the symmetric case is treated similarly). Then $k>l$, and if we write
$$
w_k(V^n)=v^n_k+\De_k,\quad w_l(V^n)=v^n_l+\De_l,
$$
our statement follows from the inequality $\De_k\geq \De_l$
(see Lemma~2) and from the equalities $v^{n+1}_k=w_k(V^n)$ and 
$v^{n+1}_l=w_l(V^n)$ (see inequalities (\ref{1.4.1})). \qed \\

Thus, we get a nondecreasing sequence of subsets of
$[1,N]\times [1,N]$:
$$
{\cal N}(V^0) \subset {\cal N}(V^1) \subset \dots \subset {\cal N}(V^n) 
\subset \cdots
$$
Since the set $[1,N]\times [1,N]$ is finite, there exists a $n_2$ and a subset ${\cal N}(V)^*$ of $[1,N]\times [1,N]$ such that
\begin{equation}
\label{103}
\nonumber
{\cal N}(V^n)={\cal N}(V)^*,\quad n\geq n_2.
\end{equation}

We again assume that $n_2=0$ and consider the set
$$
{\cal M}=[1,N]\times [1,N]\setminus {\cal N}(V)^*.
$$

By construction, this set has the following property: for any $n\geq 0$, $v^n_k$ and $v^n_l$ influence each other if and only if
$$
(k,l)\in {\cal M}.
$$

Hence,
\begin{equation}
\label{1.4.3}
J(v_k^n)=\{l\in[1,N]:\;(k,l)\in{\cal M}\},\quad k\in[1,N], \quad n\geq 0.
\end{equation}

Note that the set $J(v^n_k)$ does not depend on $n$; 
denote it $J(k)$ and let $I(k)$ be the cardinality of $J(k)$.

It is clear that, for any $k\in[1,N]$, the set $J(k)$ has the form 
$[k-\mu(k),k+\nu(k)]$, where $\mu(k),\nu(k)\geq 0$ and $\nu(k)+\mu(k)+1=I(k)$.

Introduce an $N\times N$ matrix $T$ as follows: $t_{k,l}=1/I(k)$ if $(k,l)\in {\cal M}$ and $t_{k,l}=0$ otherwise.

It follows from (\ref{1.4.3}) that
$$
\Phi(V)=(E_N+hT)V,
$$
where $E_N$ is the unit $N\times N$ matrix.

Hence,
\begin{equation}
\label{1.4.4}
V^{n}=(E_N+hT)^n V^0,\quad n\geq 0.
\end{equation}

Let us show that the spectrum of the matrix $T$ is real.

Represent $T=SU$, where $S$ is a diagonal matrix with positive diagonal elements,
$$
S=\mbox{diag}\left(\frac{1}{I(1)},\dots,\frac{1}{I(N)}\right),
$$
and entries $u_{k,l}$ of $U$ are as follows:
$u_{k,l}=1$ if $(k,l)\in {\cal M}$ and $u_{k,l}=0$ otherwise.
Clearly, $U$ is symmetric.

Then, $T$ is conjugate to
$$
S^{-1/2}TS^{1/2}=S^{-1/2}SUS^{1/2}=S^{1/2}US^{1/2},
$$
but the last matrix is symmetric:
$$
(S^{1/2}US^{1/2})^*=(S^{1/2})^*U^*(S^{1/2})^*=S^{1/2}US^{1/2}.
$$

Hence, the spectrum of $T=SU$ (and so the spectrum of $E_N+hT$) is real.

The $k$th row of the matrix $T$ has the form
$$
\left(0,\dots,0, \frac{1}{I(k)},\dots,\frac{1}{I(k)},0,\dots,0\right),
$$
where the number of nozero entries is precisely $I(k)$.

This means that $T$ is stochastic. A classical result states that $T$ has an eigenvalue 1 and all other eigenvalues $\lambda$ satisfy the inequality $|\lambda|\leq 1$.

Hence, the eigenvalues of $T$ are real and belong to $[-1,1]$, which implies that if $h\in(0,1)$, then the eigenvalues of $E_N+hT$ are positive.

In this case, any bounded sequence $V^n$ that satisfies (\ref{1.4.4}) tends to a vector $W$ such that $W=(E_N+hT)W$.  To show this,
consider a Jordan form $J$ of the matrix $E_N+hT$:
$$
J=\mbox{diag}(J_1,\dots,J_l),
$$
where $J_1,\dots,J_l$ are Jordan blocks.

Let us assume that $J_1$ is a $d\times d$ block corresponding to
an eigenvalue $\lambda$ and $d>1$ (the case $d=1$ is trivial), i.e.,
$$
J_1=\left(
\begin{array}{ccccc}
\lambda&1&0&\ldots&0\\
0&\lambda&1&\ldots&0\\
\vdots&\vdots&\vdots&\ddots&\vdots\\
0&0&0&\ldots&1\\
0&0&0&\ldots&\lambda
\end{array}
\right).
$$

Let ${k\choose j}$ be the binomial coefficients,
$$
{k\choose j}=\frac{k!}{j!(k-j)!}.
$$

If $k\geq d-1$, then
$$
J_l^k=
\left(
\begin{array}{ccccc}
\lambda^k&k\lambda^{k-1}&\frac{k(k-1)}{2}\lambda^{k-2}&\ldots&{k\choose d-1}\lambda^{k-d+1}\\
0&\lambda^k&k\lambda^{k-1}&\ldots&{k\choose d-2}\lambda^{k-d+2}\\
\vdots&\vdots&\ddots&\vdots&\vdots\\
0&0&0&\ldots&\lambda^k
\end{array}
\right).
$$

Hence,
$$
v^k_1=\lambda^k v^0_1+k\lambda^{k-1} v^0_2+\dots,
$$
$$
v^k_2=\lambda^k v^0_2+\dots,\quad\dots,\quad v^k_d=\lambda^k v^0_d.
$$

It follows that if the sequence $V^n$ is
bounded and $\lambda>1$, then $v^0_1=\dots= v^0_d=0$.

If $\lambda=1$, then $v^0_2=\dots= v^0_d=0$, and if we denote
$u=(v^0_1,0,\dots,0)$, then $J_1^n u=u$ for all $n\geq 0$.

Finally, if  $\lambda<1$, then $v^n_1,\dots, v^n_d\to 0$ as $n\to\infty$.

Of course, similar statements hold for all Jordan blocks.

This implies that if the sequence $V^n$ is
bounded, then we can represent $V^0$ in the form $W_1+W_2$ such that
$(E_N+hT)W_1=W_1$ and $(E_N+hT)^nW_2\to 0$ as $n\to\infty$.
This completes the proof.  \qed \\
\bigskip

{\bf 4. Numerical example}
\bigskip

 A numerical simulation 
 shows that the final outcome of an election process may change 
 depending on the level of interaction of the society. This has the interesting interpretation that a society is a complex entity which cannot be reduced to the simple union of many individuals: beliefs in the society evolve differently depending on the quality and level of mutual influence, which in turn is highly dependent on 
 technology and on the possible existence of rules that limit the circulation of information.

In this example, we take $N=100$, $h=0.1$, and
the initial profile of opinions is as follows:

$$v_k=-0.6, \quad 0\leq k <20;$$
$$v_k=-0.4, \quad 20\leq k <48;$$
$$v_k=-0.01, \quad 48\leq k <60;$$
$$v_k=0.1, \quad 60\leq k <90;$$
$$v_k=0.2, \quad 90\leq k \leq 100.$$
\medskip

First we take $\ep=0.45$ (high level of interaction). Numerical simulation shows
that at step 27, a clustering equilibrium is
reached, where opinion 1 achieves majority (at the equilibrium,
$v_k=-1$ for $k=1,\dots,47$ and $v_k=1$ for $k=48,\dots,100$).
\medskip

After that, we take $\ep=0.05$ (low level of interaction)
and the same initial profile. At step 49, a clustering equilibrium is
reached, where the opposite opinion $-1$ achieves majority
(at the equilibrium,
$v_k=-1$ for $k=1,\dots,59$ and $v_k=1$ for $k=60,\dots,100$).

\bigskip 

{\bf References}
\bigskip

1. M.L. Bertotti and M. Delitala,
 Cluster formation in opinion dynamics: a qualitative analysis,
 {\em Z. Angew. Math. Phys.}, {\bf 61}, 583--602, 2010.
 
 2. D. Borra and T. Lorenzi,
 Asymptotic analysis of continuous opinion dynamics models under bounded confidence,
 {\em Commun. Pure Appl. Anal.}, {\bf 12}, 1487--1499, 2013.
 
 3. F. Ceragioli and P. Frasca,
  Continuous and discontinuous opinion dynamics with bounded confidence,
  {\em Nonlinear Anal. Real World Appl.},  {\bf 13}, 1239--1251, 2012.
  
  4. S. Chatterjee,
 Reaching a consensus: Some limit theorems,
 {\em Proc. Int. Statist. Inst.}, 159--164, 1975.
 
  5. S. Chatterjee and E. Seneta,
 Toward consensus: some convergence theorems on repeated averaging,
  {\em J. Appl. Prob.}, {\bf 14}, 89--97, 1977.
 
6. J.C. Dittmer,
  Consensus formation under bounded confidence,
 {\em Nonlinear Analysis}, {\bf 47}, 4615--4621,  2001.
 
 7.  S. Fortunato,
 The Krause -- Hegselmann consensus model with discrete opinions,
 {\em Internat. J. Modern Phys. C}, {\bf 15}, 1021--1029, 2004.
  
  8. J.R.P. French,
  A formal theory of social power,
  {\em Psychological Review}, {\bf 63}, 181--194, 1956.
 
 9. G. Deffuant, D. Neau, F. Amblard, and G. Weisbuch,
  Mixing beliefs among interacting agents,
 {\em Advances in Complex Systems}, {\bf 3}, 87--98, 2000.
  
  10. G. Deffuant, G. Weisbuch,  F. Amblard, and G.P. Nadal,
 Interacting agents and continuous opinion dynamics,
 in: {\em Heterogenous Agents, Interactions and Economic Performance,} 
 Lecture Notes in Economics and Mathematical Systems, {\bf  521},
 Springer, Berlin,  2003, pp. 225--242.
 
 11. M.H. De Groot,
 Reaching a consensus,
{\em J. Amer. Statist. Assoc.}, {\bf 69}, 118--121, 1974.
  
  12. F. Harary,
 A criterion for unanimity in French's theory of social power,
 in: {\em Studies in Social Power}, Institute for Social Research, Ann Arbor,
 1959, pp. 168--182.
 
 13. P. Hegarty and E. Wedin,
  The Hegselmann -- Krause dynamics for equally spaced agents,
 {\em J. Difference Equ. Appl.}, {\bf 22}, 1621--1645, 2016.
  
  14. R. Hegselmann and A. Flache,
 Understanding complex social dynamics: a plea for cellularautomata based modelling,
  {\em Journal of Artificial Societies and Social Simulation}, {\bf 1}, 1--1, 1998.
  
  15. R. Hegselmann and A. Flache,
  Opinion dynamics and bounded confidence: Models, analysis and simulation,
  {\em Journal of Artificial Societies and Social Simulation}, {\bf 5}, 1--33, 2002.
  
  16. R. Hegselmann and U. Krause,
  Opinion dynamics under the influence of radical groups, charismatic leaders, and other constant signals: a symple unifying model,
 {\em Networks and Heterogeneous Media}, {\bf 10}, 477--509, 2015.
 
 17. J. Hajnal, J. Cohen, and C.M. Newman,
 Approaching consensus can be delicate when positions harden,
 {\em Stochastic Proc. and Appl.}, {\bf 22}, 315--322, 1986.
  
  18. P.E. Jabin and S. Motsch,
  Clustering and asymptotic behavior in opinion formation,
  {\em J. Differential Equations}, {\bf 257}, 4165--4187, 2014.
  
  19. U. Krause,
A discrete nonlinear and nonautonomous model of consensus formation,
 in: {\em Communications in Difference Equations},
 Gordon and Breach Publ., Amsterdam, 1997, pp. 227--236.
 
 20. U. Krause,
Soziale Dynamiken mit vielen Interakteuren. Eine Problemskizze,
in: {\em Modellierung und Simulation von Dynamiken mit vielen interagierenden Akteuren}, Universitat Bremen, 1997, pp. 37--51.

21. U. Krause,
 Compromise, consensus, and the iteration of means,
 {\em Elem. Math.}, {\bf 64}, 1--8, 2009.
 
 22. U. Krause,
 {\em Positive Dynamical Systems in Discrete Time. Theory, Models, and Applications},
De Gruyter, Berlin, 2015.
 
  23. S. Kurz and J. Rambau,
  On the Hegselmann -- Krause conjecture in opinion dynamics,
  {\em Journal  Difference Equ. Appl.}, {\bf 17}, 859--876, 2011.
  
  24. K. Lehrer,
  Social consensus and rational agnoiology,
  {\em Synthese}, {\bf 31}, 141--160, 1975.
  
  25. K. Lehrer and C.G. Wagner,
  {\em Rational Consensus in Science and Society},
  Dordrecht: D. Reidel Publ. Co., 1981.
 
  26. J. Lorenz,
  A stabilization theorem for dynamics of continuous opinions,
  {\em Physica A}, {\bf 355}, 217--223, 2005.
  
  27. J. Lorenz,
  Continuous opinion dynamics under bounded confidence: a survey,
  {\em International Journal of Modern Physics C}, {\bf 18}, 1819--1838, 2007.
 
 28. W. Ren and Y. Cao,
 {\em Distributed Coordination of Multi-agent Networks. 
 Emergent Problems, Models, and Issues}, Springer,  2011.
 
 29. C.G. Wagner,
Consensus through respect: a model of rational group decision-making,
{\em Philosophical Studies}, {\bf 34}, 335--349, 1978.
    
  30. E. Wedin and P. Hegarty,
The Hegselmann -- Krause dynamics for the continuous-agent model and a regular opinion function do not always lead to consensus,
 {\em IEEE Trans. Automat. Control}, {\bf 60}, 2416--2421, 2015.
 
  31. S. Wongkaew, M. Caponigro and A. Borzi,
  On the control through leadership of the Hegselmann -- Krause 
  opinion formation model,
{\em Mathematical Models and Methods in Applied Sciences}, {\bf 25}, 565--585, 2015.

\end{document}